\documentclass[12pt,a4paper]{article}
\usepackage{latexsym,epsfig,fullpage,amsfonts,amssymb,amsmath,
amscd,epic,graphics,
theorem}
\theoremstyle{plain}
\newtheorem{lemma}{Lemma}

\newtheorem{remark}[lemma]{Remark}

\newtheorem{definition}[lemma]{Definition}

{\theorembodyfont{\rmfamily} 
  
\begin{document}
\baselineskip=15pt
\newcommand{\pperp}{\hbox{$\perp\hskip-6pt\perp$}}
\newcommand{\ssim}{\hbox{$\hskip-2pt\sim$}}
\newcommand{\N}{{\mathbb N}}\newcommand{\Delp}{{\Pi}}
\newcommand{\A}{{\mathbb A}}
\newcommand{\Z}{{\mathbb Z}}
\newcommand{\R}{{\mathbb R}}
\newcommand{\C}{{\mathbb C}}
\newcommand{\Q}{{\mathbb Q}}
\newcommand{\PP}{{\mathbb P}}
\newcommand{\mnote}{\marginpar}
\newcommand{\Proj}{{\mathbf{Proj}}}
\newcommand{\oeps}{{\overline\eps}}
\newcommand{\oDel}{{\widetilde\Del}}
\newcommand{\real}{{\operatorname{Re}}}
\newcommand{\conv}{{\operatorname{conv}}}
\newcommand{\Span}{{\operatorname{Span}}}
\newcommand{\Ker}{{\operatorname{Ker}}}
\newcommand{\Hyp}{{\operatorname{Hyp}}}
\newcommand{\Fix}{{\operatorname{Fix}}}
\newcommand{\sign}{{\operatorname{sign}}}
\newcommand{\Tors}{{\operatorname{Tors}}}
\newcommand{\oi}{{\overline i}}
\newcommand{\oj}{{\overline j}}
\newcommand{\ob}{{\overline b}}
\newcommand{\os}{{\overline s}}
\newcommand{\oa}{{\overline a}}
\newcommand{\oy}{{\overline y}}
\newcommand{\ow}{{\overline w}}
\newcommand{\ou}{{\overline u}}
\newcommand{\ot}{{\overline t}}
\newcommand{\oz}{{\overline z}}
\newcommand{\bw}{{\boldsymbol w}}
\newcommand{\bx}{{\boldsymbol x}}
\newcommand{\bu}{{\boldsymbol u}}
\newcommand{\bz}{{\boldsymbol z}}
\newcommand{\eps}{{\varepsilon}}
\newcommand{\proofend}{\hfill$\Box$\bigskip}
\newcommand{\Int}{{\operatorname{Int}}}
\newcommand{\pr}{{\operatorname{Pr}}}
\newcommand{\grad}{{\operatorname{grad}}}
\newcommand{\rk}{{\operatorname{rk}}}
\newcommand{\im}{{\operatorname{Im}}}
\newcommand{\sk}{{\operatorname{sk}}}
\newcommand{\const}{{\operatorname{const}}}
\newcommand{\Sing}{{\operatorname{Sing}}}
\newcommand{\conj}{{\operatorname{Conj}}}
\newcommand{\Pic}{{\operatorname{Pic}}}
\newcommand{\Crit}{{\operatorname{Crit}}}
\newcommand{\Ch}{{\operatorname{Ch}}}
\newcommand{\discr}{{\operatorname{discr}}}
\newcommand{\Tor}{{\operatorname{Tor}}}
\newcommand{\Conj}{{\operatorname{Conj}}}
\newcommand{\val}{{\operatorname{val}}}
\newcommand{\Val}{{\operatorname{Val}}}\newcommand{\GW}{{\operatorname{GW}}}
\newcommand{\defect}{{\operatorname{def}}}
\newcommand{\tmu}{{\C\mu}}
\newcommand{\ov}{{\overline v}}
\newcommand{\ox}{{\overline{x}}}
\newcommand{\tet}{{\theta}}
\newcommand{\Del}{{\Delta}}
\newcommand{\bet}{{\beta}}
\newcommand{\kap}{{\kappa}}
\newcommand{\del}{{\delta}}
\newcommand{\sig}{{\sigma}}
\newcommand{\alp}{{\alpha}}
\newcommand{\Sig}{{\Sigma}}
\newcommand{\Gam}{{\Gamma}}
\newcommand{\gam}{{\gamma}}
\newcommand{\Lam}{{\Lambda}}
\newcommand{\lam}{{\lambda}}
\newcommand{\SC}{{SC}}
\newcommand{\MC}{{MC}}
\newcommand{\nek}{{,...,}}
\newcommand{\cim}{{c_{\mbox{\rm im}}}}
\newcommand{\mathto}{\mathop{\to}}
\newcommand{\op}{{\overline p}}

\newcommand{\w}{{\omega}}

\def\Top{\operatorname{Top}}
\def\Def{\operatorname{Def}}
\def\Im{\operatorname{Im}}
\def\Wu{\operatorname{Wu}}
\def\inj{\operatorname{inj}}

\let\+\sqcup

\title{Overview of topological properties of real algebraic surfaces}
\author{Viatcheslav Kharlamov}
\date{}
\maketitle

\vskip1cm

\begin{abstract}
These notes reproduce the content of a short, 50-minutes,
survey talk given at the Nice University in September, 2004.
We added a few topics that have not been touched on in the lecture
by lack of time.
\end{abstract}

\section{Introduction}\label{intro}

%
%
%

Topology of real algebraic varieties is a broad subject. Thus,
it is reasonable to specify the "level" of objects and the
goal of study. In what concerns the level, one may distinguish
between affine, projective and abstract varieties, and, from a certain
point of view, it is natural to start with abstract varieties, and
then descend to projective and affine ones. As to the goal, I'll give
preference to those "real results" that require "complex proofs" (even
though in what follows the proofs will almost always stay behind the
scenes) and, moreover, admit "complex statements".

Another major, and traditional, simplification is to consider
{\bf nonsingular} varieties, at least at the first stage. Certainly,
a complete separation from the singular world is never
possible. However, de facto, in all the cases when a complete
understanding was achieved in the nonsingular case, it turned out
that the singular case could be treated, at least in principle, by
similar methods.

In such a setting, it is natural to consider not only algebraic, but
arbitrary K\"ahler compact complex manifolds, and to call a complex
manifold $X$ {\it real}, if it can be equipped with a {\it real
structure}, that is an anti-holomorphic involution $c:X\to X$.
Real points are then, by definition, the fixed points of the real
structure. We denote by $X_\R$ the set of real points, or the {\it real
part}, of $X$. For the sake of symmetry, we often denote $X$ by $X_\C$.

Of course, the principal source of examples is given by nonsingular
varieties defined by systems of real polynomial equations; in these
$c$ is the complex conjugation. Note that by our convention a real
variety is nonsingular if it does not have singular points, be they real
or imaginary.

We consider real algebraic varieties up to the following equivalence
relations: diffeomorphism of real part, diffeomorphism of real
structure, and deformation of variety together with real structure.
As usual, by an {\it elementary real deformation} of a {
real variety we mean a smooth, differentially locally trivial, family of
{
real varieties (say, an equivariant deformation in the sense of
Kodaira-Spencer). Two real varieties are called
real deformation equivalent if
there exists a chain of elementary deformations connecting them.
{\it Topology of the real part and that of the real structure are
preserved under deformation.} This phenomenon is one of the main
motivations for the study of topological properties of real varieties.
Note in advance that in many important cases the topology of real structure
(which includes, in fact, the topology of real part) determines the
deformation class.

A fundamental example of deformation is a small variation of a nonsingular
system of polynomial equations (that is, a system whose Jacobian has maximal
rank
at each solution of the system) or, more generally, any variation which is represented in
the total space of systems of a given number of equations (equal to the
codimension of variety) of a given degree in a given
number of variables by a smooth path in the complement of the discriminant
locus, that is a smooth path avoiding singular systems.
However, other deformations may exist as well; thus, surfaces of degree
$\ge 5$ in $\PP^3$ have big deformations failing to be embedded in $\PP^3$
while surfaces of degree $4$ have small K\"ahler deformations failing to be
projective. By contrast, any deformation of a surface of degree $3$ is
realized in $\PP^3$.

Our choice of surfaces as the topic of this talk is motivated by the fact
that this is the first nontrivial case with respect to the above equivalence
relations. In fact, in dimension zero the topology is determined by the natural
number $b_0(X_\R)$ satisfying the following relations
\begin{equation}\label{smip}
b_0(X_\R)\le b_0(X_\C),\quad b_0(X_\R)=b_0(X_\C)\mod
2 
\end{equation}
(here and in what follows $b_i(\cdot)$ denotes the rank of $H_i(\cdot\, ;\Z/2\Z)$,
so that $b_0$ is nothing but the number of connected components).
Similarly, in dimension one everything is determined by the number $b_0(X_\R)$
(or, equivalently, by $b_1(X_\R)=b_0(X_\R)$), and the only relations linking
this number with the invariants of $X_\C$ are as follows:
\begin{equation}\label{smic}
b_0(X_\R)\le g(X_\C)+1,\quad \text{and}\quad b_0(X_\R)=g(X_\C)+1\,\, \text{\rm
mod}\,
2 \,\,\text{\it if}\,\, X_\R\,\text{\it divides}\, X_\C.
\end{equation}
Here, $g$ is the usual genus, $g(X_\C)=\frac12 b_1(X_\C)$ if $X$
is irreducible; otherwise, $g(X_\C)+1$ is the sum of $g+1$ over all
irreducible components. The condition that $X_\R$ divides $X_\C$ is
equivalent to the orientability of the quotient $X_\C/c$.
The above relations demonstrate a general phenomenon: $\Top (X_\R)\le
\Top (X_\C)$ (the topology of a real variety is bounded by that of its
complexification).

Moreover, in dimensions less than two topology, and even homology
determines the deformation equivalence classes: $\Top (X_\C)=\Def (X_\C)$
over $\C$ and $\Top (X_\C,c)=\Top (X_\C,X_\R)=\Def (X_\C,c)$
over $\R$. For example, two real irreducible curves are deformation
equivalent, if and only if they have the same number of real components,
$b_0(X_\R)$, and they both either divide $X_\C$ or not (in other words,
either, in both cases, $X_\R$ is homologous to $0$ in $X_\C$, or not).

An advantage of dimension two is that algebraic topology still suffices
to determine $\Top X_\R$, which is no more the case in dimensions $\ge 3$
(that is why the dimensions $\ge 3$ are still very far from being well
understood; cf., however, J.~Koll\'ar's papers \cite{Ko1, Ko2} and
references in \cite{Bour}
for some nontrivial partial results).

Many of the tools used in dimension $\le 2$ can be extended to higher
dimensions, and when consideration of higher dimensions does not lead to
complications we present our results in full generality.

For lack of time, we do not discuss arrangements of curves on surfaces,
construction of surfaces (and curves on surfaces) with prescribed topology,
or enumerative results. We also will not have time to discuss history of
the subject, but it is worth mentioning that many of the results presented
below emerged due to collective efforts of many insiders (in particular,
a group of Russian mathematicians inspired by I.~G.~Petrovsky and
V.~I.~Arnol'd in Moscow, D.~A.~Gudkov in the former Gor'kii, and
V.~A.~Rokhlin in the former Leningrad; the author gratefully dedicates
these notes to the memory of the latter). A reader interested to know
better who, when and how discovered these results is invited to consult
\cite{DK} and the references therein.

\section{Topology of real varieties versus topology of their complexifications}

\subsection{ Smith theory}

Smith theory provides the following relations valid for all
dimensions:
\begin{equation}\label{smith}
\sum b_i(X_\R)\le \sum b_i(X_\C),\quad \sum b_i(X_\R)=\sum b_i(X_\C)\mod
2, 
\end{equation}
where $b_i(\cdot)=\dim H_i(\cdot\, ;\Z/2\Z)$.
Behind these relations there are such
useful tools as Smith's exact sequence and Kalinin's spectral
sequence, see, for example, \cite{DK}. The latter starts with
$$
E^1_*=H_*(X_\C;\Z/2\,\Z), \quad E^2_*=\Ker(1+c_*)/\Im(1+c_*),
$$
and
converges to $H_*(X_\R;\Z/2\,\Z)$ (here and in what follows $c_*$ states for the homomorphism
$H_*(X_\C;\Z/2\,\Z)\to H_*(X_\C;\Z/2\,\Z)$ induced by the real structure
$c:X_\C\to X_\C$). Already the existence of such a
spectral sequence implies (\ref{smith}). According to the above
formula for $E^2_*$, there is a "stronger" inequality
\begin{equation}\label{swan}
\sum b_i(X_\R)\le \dim H_1(\Z/2\,\Z, H_*(X_\C;\Z/2\,\Z)),
\end{equation}
where $H_1(\Z/2\,\Z, H_*(X_\C;\Z/2\,\Z))=\Ker(1+c_*)/\Im(1+c_*)$.

There are two important classes of real varieties enjoying special features,
viz.
\begin{itemize}
\item
$M$-{\it varieties}, i.e. varieties for which $E^1_*=H_*(X_\R;\Z/2\,\Z)$
(which is equivalent to $\sum b_i(X_\R)= \sum b_i(X_\C)$, the extremal case
of (\ref{smith}), or to the existence of $c$-invariant cycles in any homology
class in
$H_*(X_\C;\Z/2\,\Z)$); here $M$ stands for "maximal";
\item
and
$GM$-{\it varieties}, i.e. varieties for which $E^2_*=H_*(X_\R;\Z/2\,\Z)$ (the
extremal case of (\ref{swan}), which is equivalent to
the existence of $c$-invariant cycles in any of the classes in $H_*(X_\C;\Z/2\,\Z)$
fixed by $c_*$; note that all real surfaces with $\pi_1(X_\C)=1$ are $GM$); here
$GM$ stands for "Galois maximal".
\end{itemize}

In dimensions $0$ and $1$ the relations (\ref{smith})\,  yield (\ref{smip})\, and
(\ref{smic}), except for the congruence in (\ref{smic}) which is a {\rm mod} $4$
relation for $\sum b_i$. The latter can be generalized in the
following way: {\it if $X$ is a $GM$-variety of odd dimension, $\dim
X=2k+1$,
and $X_\R$ is homologous to $0$ in $X_\C$, then}
$$
\frac12\sum b_i(X_\R)=\frac12\sum b_i(X_\C)\mod 2.
$$
Note that both parts of the above congruence are integers. Since
$X_\R\sim 0$, for any $x\in H_{2k+1}(X_\C;\Z/2)$ one has
$$(x,c_*x)=(x,[X_\R])=0\mod 2,
$$
which implies that, besides the summands with constant action of
$c$,
any irreducible orthogonal decomposition of
$c_*:H_*(X_\C;\Z/2)\to H_*(X_\C;\Z/2)$ contains
only irreducible components
of rank $4$ (each such component is a permutation of two $\Z/2\,\Z$-planes). Thus,

$$
\sum b_i(X_\R)=
\dim \Ker(1+c_*)/\Im(1+c_*)=
\sum b_i(X_\C) - 4p,
$$
where $p$ is the number of the above irreducible components of
rank $4$.

Smith theory helps to answer the following questions: {\it Is $X_\R$ non empty
(existence of real solutions)? How many connected components
does it have?}  and, more generally, {\it What are the possible values of
$b_i(X_\R)$?}

Already the congruence part of (\ref{smith})\, gives
a simple, often useful, sufficient condition: {\it $X_\R$ is non-empty, if
$\sum b_*(X_\C)$ is odd}. Note that this condition does not depend on the choice
of real structure.

Smith theory is not well adapted to work with individual Betti numbers, and so
in dimensions $>2$ it is hard to get more information and to answer the other
questions using Smith theory alone. Hopefully, in dimension $2$ it is sufficient
to perform homological calculations in Kalinin's spectral sequence
(or in Smith's exact sequence) and to add the information coming from the Lefschetz
trace formula. In particular, in the case of real surfaces with $\pi_1(X_\C)=1$
we get the following formulae:
$$
b_0(X_\R)=\frac12(\sum
b_i(X_\C)
-b_2^{-1}-a),
$$

$$
b_1(X_\R)=
b_2-b_2^{+1}-a,
$$
{\it where $b_2^{\pm 1}$ are the dimensions of the eighenspaces of the involution}
$$
c_*:H_2(X_\C;\R)\to
H_2(X_\C;\R),
$$
and $a$ is the number of nontrivial components
in an irreducible $\Z/2\Z$-vector space decomposition of $c_*:H_*(X_\C;\Z/2)\to H_*(X_\C;\Z/2)$
(each such component is generated by two elements permuted by $c_*$).

We recall that the above definition of $GM$-variety is equivalent to
the existence of an equivariant cycle in each invariant
$\Z/2\,\Z$-homology class, and the definition of $M$-variety
is equivalent to the existence of an equivariant cycle in each
$\Z/2\,\Z$-homology class. The latter happens, for example, if all
$\Z/2\,\Z$-homology classes of $X_\C$ are algebraic and have a real
representative. This is the case for projective spaces and Grassmann
varieties equipped with their tautological real structures. Many
other special varieties also have this maximality property. An
important class of varieties for which the maximality question is
open is provided by the discriminants of polynomials in three or
more (homogeneous) variables. For three variables, this reduces to the
question of maximality of the space of all singular plane curves of
a given degree. In degrees $1$, $2$, and $3$ discriminant is indeed an
$M$-variety, as can be verified by a more or less straightforward
calculation using
the Alexander-Pontryagin duality (one should first compute the Betti
numbers of the space of nonsingular curves). But already in degree $4$,
the question is open.

The relations (\ref{smith})\, apply, in fact, to any
finite-dimensional space with involution, and, in particular,
to singular varieties. In many cases, for example in the case
of projective hypersurfaces, the corresponding upper bounds are
the best known one. These relations can also be applied to
semi-algebraic sets; to do this, it suffices to replace such a set
by its tubular neighborhood and then apply (\ref{smith}) to its boundary,
which is a hypersurface. Y.~Laszlo and C.~Viterbo \cite{LV} recently
addressed the following question: how to bound the total Betti number
of a nonsingular real projective variety $X$ in terms of its degree $d$
and dimension $n$. Combining (\ref{smith})\,
with some inequalities due to Demailly-Peternell-Schneider they proved that
\begin{equation}
\sum b_i(X_\R)\le 2^{n^2+2} d^{n+1},
\end{equation}
while all the previously known estimates (like those of R.~Thom
\cite{Th} and J.~Milnor \cite{Mi}) were of the type
\begin{equation}
\sum b_i(X_\R)\le C d^{2n+1}.
\end{equation}
Using Lefschetz pencil, resp. Morse function type arguments
applied to $X_\C$, resp. $X_\R$, one can improve the leading coefficient
in the Laszlo-Viterbo bound and
get
\begin{equation}\label{coef}
\sum b_i(X_\R)\le d^{n+1} + O(d^n),\quad  \sum b_i(X_\C)\le d^{n+1} +
O(d^n)
\end{equation}
(without appealing to Demailly-Peternell-Schneider
inequalities).
What is an optimal choice of $O(d^n)$? I do not know. At least,
\begin{equation}\label{coeff}
\sum b_i(X_\R)\le \sum b_i(X_\C)\le d^{n+1}
\quad\text{if}\quad d>1,
\end{equation}
and it seems reasonable to expect that $d^{n+1} + O(d^n)$ in
(\ref{coef}) can be replaced by the polynomial
$d^{n+1}+\sum_{k=0}^n a_kd^k$
representing
the total Betti number of a degree $d$ nonsingular hypersurface in
$\C\Bbb P^{n+1}$ (cf. (\ref{total}) for the case $n=2$).
The key point in
the proof of (\ref{coef}) and (\ref{coeff})
is a similar bound,
\begin{equation}\label{dual}
d^*\le d(d-1)^n=d^{n+1}+O(d^n),
\end{equation}
for the degree $d^*$
of the
variety $X^*$ projectively dual to $X$ (I am grateful to F.~Zak who
explained me how such a general bound for $d^*$
is deduced from the, classical, computation of $d^*$ in the special
case of hypersurfaces; he also proposed to replace $O(d^n)$ by $0$
in (\ref{coef})). Indeed, the number of singular fibers of
the pencil of hyperplane sections of $X_\C$ and, respectively, the number of
singular values of the linear Morse function on an affine part of $X_\R$
are bounded by
$d^*$, so that an induction on the dimension $n$ of $X_\R$ gives
a sequence $x_n=\sum b_i(X_\R)$ (resp. $y_n= \sum b_i(X_\C)$)
with the property
$$ x_{n}-2x_{n-1}+x_{n-2}\le  d^*
\quad(\text{respectively}\quad
y_{n}-2y_{n-1}+y_{n-2}\le
d^*)
$$
(here $X=X_n, X_{n-1},\dots, X_0$ is the sequence formed by $X$ and its
consecutive hyperplane sections). Combined with (\ref{dual}), this
yields (\ref{coef}). For $d>2$ and $n>2$ the bound (\ref{coeff}) follows
from the inequalities $y_{n}\le d^*+2y_{n-1}\le d(d-1)^{n}+2d^{n}\le
d^{n+1}$, and if $d=2$ or $n=1,2$, then  it is easy to prove the bound by
an ad hoc argument.

Returning to the Smith bound, let us forewarn that it is impossible
in general
to replace $b_i(\cdot)=\dim H_*(\cdot\,;\Z/2\,\Z)$ in the inequality
(\ref{smith}) by the ordinary Betti numbers $\beta_i(\cdot)=\dim H_*(\cdot\,;\mathbb
Q)$. For example, there exist real Enriques surfaces $X$ with the
real part consisting of two real components, one homeomorphic to a
torus and another to a connected sum of $10$ real projective planes,
while for such real surface $ \dim H_*(X_\C;\mathbb Q)=12<14=
\dim H_*(X_\R;\mathbb Q)$.

\subsection{ Higher order congruences}\label{cong}

Higher order congruences can be found basing on Smith theory and
arithmetic of integral quadratic forms. Here is a typical example:
{\it if $\dim X=2k$ and $X_\R$ is $\Z/2$-homologous to the
middle dimensional Wu class of $X_\C$, then}
$$
\chi (X_\R)=\sigma (X_\C)\mod 8
$$
(here $\chi$ is the Euler characteristic and $\sigma$ is the
signature).
Its one-line proof given below is a model for finding other
higher order congruences. It is based on the Lagrangian property
of the real part and the Lefschetz-Hirzebruch signature formula:
$$
\chi(X_\R)=(-1)^k\sigma(c)=_8 (-1)^k(\Wu_c, c\Wu_c)=
$$
$$
(-1)^k(\Wu_X,
c\Wu_X)=(\Wu_X,\Wu_X)=_8\sigma(X_\C)
$$
(here $\Wu_c$ is the Wu integral characteristic class of the quadratic
form $(x,cy)$ on $H_{2k}(X_\C,\Z)/\text{Tors}$,
$\sigma(c)$ is its signature,
$\Wu_X$ is an integral algebraic representative of the $\Wu$-class of
$X_\C$, and $=_8$ stands for congruence modulo $8$). It should be mentioned that in the case when $X_\C$ is a
complete intersection there are various methods for explicit
computation of $\sigma(X_\C)$ (see, for example,
\cite{DK} and references therein).

There is a series of congruences refining the above one (see, for
example, \cite{DK} and references therein); the two simplest of them
are as follows:

{\it if $\dim X=2k$ and $X_\R$ is a $M$-variety, then}
$$
\chi (X_\R)=\sigma (X_\C)\mod 16;
$$

{\it if $\dim X=2k$ and $X_\R$ is a $(M-1)$-variety, then}
$$
\chi (X_\R)=\sigma (X_\C)\pm 2\mod 16.
$$

Note that $\sigma(X_\C)=(-1)^k\chi (X_\C)\mod4$
(which most easily can be seen from the Hodge decomposition),
which yields an analog of congruence (\ref{smic}) in even dimensions:

$\chi (X_\R)=(-1)^k\chi (X_\C)\mod 4$ {\it if $X_\R$ is
$\Z/2$-homologous to the middle dimensional Wu class of $X_\C$}.

\subsection{An application of Hodge theory and some other inequalities}

As is shown in \cite{Kh}, from the Hodge decomposition and the Lefschetz
formula it follows that
\begin{equation}\label{kh}
\vert\chi(X_\R)-1\vert\le h^{k,k}(X_\C)-1
\quad \text {\it if}\quad \dim X=2k
\end{equation}
(where $h^{k,k}$ is the Hodge number of bidegree $(k,k)$; various
explicit
computations of the Hodge numbers are found in \cite{DK}; an expression
for
$h^{1,1}$ in the case of surfaces in $3$-space is given below in (\ref{h11})).
See \cite{Kh2,Kh3, DK} and references therein for an odd-dimensional version
of this Comessatti-Petrovskii type inequality and
for generalizations to varieties with singularities (naturally,
in the singular case pure Hodge structure is to be replaced by mixed one).
It would certainly be nice to find other applications
of Hodge theory giving more detailed information than (\ref{kh}).
Especially challenging is to somehow relate Hodge theory with Smith theory.

Let me indicate here only a very special amusing
application of (\ref{kh}) to the case of odd dimension.
It concerns plane curves, and, more specifically, line arrangements.
We consider
a generic configuration of $2k$ real lines in the projective
plane. The number of connected components, called cells, of the complement
of the arrangement is equal then to $2k^2-k+1$.
Since the number of lines is real, the cells can be chess-board colored,
and an application of (\ref{kh}) shows that
an upper bound for the number of projective cells of one color is
$\frac32 k(k-1)+1$, so that
a bound from below for the other color is given by
$\frac12 k(k+1)$.

More special inequalities, not directly related to Hodge theory,
can be obtained using the Lagrangian property of $X_\R$. Thus, in the case of
surfaces one can easily show that

{\it the number $p_-$ of orientable
components of $X_\R$
with $\chi <0$ has an upper bound
$$
p_-\le \frac12 (\sigma^+(X_\C)-1)
$$
{\rm (}where $\sigma^+$ is the positive index of inertia of the
intersection form}{\rm ).}

Note also that in the case of surfaces the inequality (\ref{kh})
can be extended to non-K\"ahler surfaces in the form
$$
2-\dim H^2(X_\C;\R)+2h^{2,0}(X_\C)\le \chi (X_\R)\le 4+\dim
H^2(X_\C;\R)-2h^{2,0}(X_\C),
$$
which is weaker than (\ref{kh}), but differs from it only by $4$ in the
left- and right-hand parts. This difference is due to the absence of
a K\"ahler class in $H^2$ and the asymmetry $H^{1,0}=H^{0,1}-1$ in $H^1$.

It would be interesting to find analogs of (\ref{kh}) for
the signature $\sigma(X_\R)$ instead of $\chi(X_\R)$, of course
under the hypothesis that $X_\R$ is orientable. The best
bound known to me does not involve Hodge theory. It says that
\begin{equation}\label{vi}
\vert\,\sigma(X_\R)\vert \le \frac13 \,c_2^2(X_\C),
\end{equation}
and follows directly from evaluating an algebraic representative of the
second Chern class of $X_\C$ on $X_\R$, viz.
$\inj_*c_2(X_\C)\cap[X_\R]=p_1 (X_\R)\cap[X_\R]=3\sigma(X_\R).$
This bound holds under an additional hypothesis that the tangent,
or cotangent, or some other vector bundle of $X_\C$ with the same $c_2$ is
generated by its sections (this moving condition allows to put an algebraic
representative of the second Chern class in a general position with
respect to $X_\R$ and thus to get (\ref{vi})).
More general and considerably more subtle bounds for arbitrary Pontryagin
numbers can be found in a recent paper by Y.~Laszlo and C.~Viterbo \cite{LV}.

\subsection{Special surfaces}\label{types}
The above tools allow to
understand thoroughly the topology of $X_\R$ for many special
types of surfaces. For example, they lead to a complete
topological classification of $X_\R$, and even of $(X_\C,c)$,
for cubic and quartic surfaces in $\PP^3$. We describe it in terms
of generators: each topological type generates a list of its Morse
simplifications, that is the topological types obtained from the
initial one
by series of Morse surgeries decreasing the total Betti number (removing
a spherical component or contracting a handle).

{\it There are
$5$ classes of nonsingular cubics generated by $\#_7\R P^2$
and $\,\R P^2\+ S^2$ {\rm (here and in what follows $\# $ stands for a
connected and $\+$ for a disjoint sum)}, and $66$ classes of nonsingular
quartics generated by three {\rm $M$-surfaces}
$\#_{10}(S^1\times S^1)\+ S^2$, $\#_{6}(S^1\times S^1)\+ 5 S^2$,
$\#_{2}(S^1\times S^1)\+ 9 S^2$, two {\rm $(M-2)$-surfaces}
$\#_{7}(S^1\times S^1)\+ 2S^2$, $\#_{3}(S^1\times S^1)\+ 6S^2$,
and a {\rm pair of tori} $2(S^1\times S^1)$.}

Surfaces in $\R P^3$ can also be studied up to different equivalence
relations, such as: {\it ambient isotopy in $\R P^3$}, {\it rigid
isotopy} (i.e., isotopy in the class of nonsingular or, more
generally, equisingular in some appropriate sense surfaces of the same
degree), and {\it rough projective equivalence} (i.e., projective
transformation and rigid isotopy). The difference between the last two
relations is due to the fact that the group $\text{\sl PGL}(4;\R)$ of
projective transformations of $\R P^3$ has two connected components.
Of course, the transformations in the component of unity transform a
surface into a rigidly isotopic one. To what extent the
classifications up to rigid isotopies and up to rough projective
equivalence are topological is an open question for surfaces of
degree $5$ and higher, cf. the discussion in \ref{qs}.

Topologically, the non-spherical component of
the real part of a nonsingular cubic
surface is embedded in $\R P^3$ as the standard $\R P^2$ with unlinked and
unknotted handles attached. Moreover, for cubic surfaces
not only the isotopy equivalence relation, but all
the other relations
mentioned above coincide with the purely topological one.

The embedding of quartic in~$\R P^3$ is also simple: it is isotopic to a
union of ellipsoids and hyperboloids with unknotted and unlinked handles.
With one exception, the components are outside each other; in the exceptional
case the real part consists of two nested spheres. In all other cases the
isotopy type of the real part~$X_\R$ of a real quartic
surface in $\R P^3$ is determined by its topological type and contractibility
or noncontractibility of $X_\R$ in $\R P^3$. It turns out that
in the case of degree~$4$ surfaces all the four classifications
(topological, isotopic, rough projective, and rigid) are different.
Note that the only difference between rough projective
equivalence and rigid isotopy is in chirality which tells
whether or not a surface is rigidly isotopic to its mirror image.
Rough projective equivalence is discussed in \cite{Nik2} and
chirality in \cite{Khchir}.

Any nonsingular degree 4 surface in $P^3$ is a {\it $K3$-surface},
that is a compact complex
surface with $\pi_1=1$ and $c_1=0$. Other examples of $K3$-surfaces
are given by double coverings of a nonsingular quadric in $P^3$ branched
in a transverse section by a quartic, by double coverings of a
nonsingular cubic in $P^3$ branched in a transverse section by a quadric,
by transverse intersections of three quadrics in $P^5$, etc.
Classification of all the real projective $K3$-surfaces up to rough projective
equivalence can also be found in \cite{Nik2}.

The methods used in the study of real $K3$-surfaces are based on the above tools,
including Hodge theory, as well as the Torelli theorem which plays a key role.
In what concerns rigid isotopies and rough projective equivalence,
using the surjectivity of the period map, one can reduce the
study of real structures to a study of arithmetic properties of integral
lattices. Similar methods can be used to study $K3$-surfaces with simple
singularities, but this problem has never been treated systematically.

Starting with degree $5$, our knowledge is much more limited.
It is not even known what are
the extremal values of the Betti numbers of nonsingular quintics.
We only know that the maximal number of connected components is
somewhere in between $23$ and $25$ and that the maximal first
$\Z/2\Z$-Betti number is either $45$ or $47$ (for the surfaces
in the same deformation class one has $\max b_1=47$). The best
known general bounds for the Betti numbers are those given by the
inequalities described in the previous sections. For a
surface in $\PP^3$ and, more generally, for a transversal complete
intersection $X$ in~$\PP^q$, the complex ingredients of these bounds
can easily be found. To wit, if
$X$ is a complete intersection in~$\PP^{q}$ of polydegree
$(m_1,\dots,m_{q-2})$ then
$$
 b_1(X)=0,\qquad b_2(X)=\chi(X)-2,
 $$
 $$
 h^{1,1}(X)=\frac12[b_2(X)-\sigma(X)]+1
 $$
 $$
 \chi(X)=\mu_{q-2}\left(\mu_1^2-\mu_2-(q+1)\mu_1+\tfrac12q(q+1)\right),
 $$
 $$
 \sigma(X)=-\tfrac13\mu_{q-2}(\mu_1^2-2\mu_2-q-1),
$$
where $\mu_i$ is the $i$-th elementary symmetric polynomial in
$(m_1,\dots,m_{q-2})$. In particular, for a surface of degree $m$
in $\PP^3$
\begin{equation}\label{total}
\sum b_i(X_\C)=\chi(X_\C)=m^3-4m^2+6m,
\end{equation}
\begin{equation}\label{h11}
h^{1,1}(X_\C)-1=(m-1)^3-\frac13 m(m-1)(m-2).
\end{equation}

As it was already noticed, the same tools can be applied to singular
objects as well (see, for example, \cite{Kh3} and \cite{Var}). For instance,
one can use them to bound the number of real double points in the following very
simple way. In the
case of surfaces there are two types of such points, viz. solitary
points and nodes (in local coordinates their equations are
$x^2+y^2+z^2=0$ and $x^2+y^2=z^2$ respectively).
One can resolve the nodes (which is differentially equivalent to
replacing a neighborhood of a node by its perturbation $x^2+y^2=z^2+\epsilon^2$)
and replace the solitary points by spheres (which means replacing
a neighborhood of a solitary point in $X_\C$ by its perturbation
$x^2+y^2+z^2=\epsilon^2$). As a result, we obtain a $4$-manifold
diffeomorphic to the minimal desingularization $\tilde X_\C$
of $X_\C$ and an
involution on it such that the fixed point set is diffeomorphic
to a disjoined sum of the minimal desingularization $\tilde X_\R$ of $X_\R$
and $S$  spheres, where $S$ is the number of solitary points
of $X_\R$. Now, applying the Smith inequality, one gets
$$
2S+\sum b_i(\tilde X_\R)\le\sum b_i(\tilde X_\C).
$$
Thus in the case of surfaces of degree $m$ in $\PP^3$

$$
2S\le 2S+\sum b_i(\tilde X_\R)\le m^3-4m^2+6m.
$$
This implies, in particular, that
$S\le \frac12(m^3-4m^2+6m)$. Using the
congruences $\mod 16$ described
in \ref{cong}, this
 can be improved to a sharp bound: {\it the number of solitary
points of a real quartic and, more generally, of any real singular $K3$-surface,
is $\le 10$.} (This may be worth comparing with the upper bound $16$
for the number of complex nodes of a complex quartic. This bound, which
is probably due to R.W.H.T.~Hudson, was extended to any singular $K3$-surface
by V.~Nikulin \cite{Nik} who used arithmetic of integral quadratic forms. As is
well known, probably since Fresnel and Kummer, real quartics with $16$ real
nodes do exist.)

Let me notice that the frontier of our knowledge of surfaces in
$\PP^3$ is similar to the frontier between special surfaces and
surfaces of general type in the Enriques-Kodaira classification of
compact complex surfaces: surfaces of degree $\ge 5$ are of
general type while surfaces of degree $4$ are
$K3$-surfaces and surfaces of degree $3$ are rational. This
gives additional motivation to turn to real structures
on complex surfaces in various Enriques-Kodaira classes.

\section{Deformation classes}

Even the above very sketchy discussion shows that
a thorough topological study of surfaces leads unavoidably to
their study up to variation of equations and then to their study
up to deformation (see Introduction for the definition; recall that
we have chosen to work with K\"ahler surfaces).

\subsection{Quasi-simplicity}\label{qs}
As is pointed in \ref{types},
two nonsingular real cubic surfaces are real deformation equivalent
if and only if their real point sets are homeomorphic.
Furthermore, the real structures of two real nonsingular cubic surfaces
are
diffeomorphic if (and only if) the real point sets of the
surfaces are homeomorphic.
This is a manifestation of what we call the quasi-simplicity
property:
a real surface $X$ is called {\it quasi-simple} if it is
real deformation
equivalent to any other real surface $X'$ such that, first, $X'$ is
deformation equivalent to $X$ as a complex surface, and, second, the
real structure of $X'$ is diffeomorphic to the real structure of $X$.

In fact, all rational real surfaces are quasi-simple.
For $\mathbb R$-minimal (i.e., minimal over $\mathbb R$) rational surfaces this result is essentially due
to Comessatti, Manin, and Iskovskikh (see e.g. the survey \cite{Man}).
In full generality this is proved in \cite{DKrat}, where, in addition,
it is shown that the real deformation type of a real rational surface is
determined by certain homological data.

Ruled $\mathbb C$-minimal surfaces of any genus are also quasi-simple, see \cite{W}.
Another class of surfaces whose real deformation theory is well understood
is formed by minimal surfaces of Kodaira dimension $0$. This class consists
of Abelian, hyperelliptic, $K3$-, and Enriques surfaces. They are all
quasi-simple (see \cite{DIK} and \cite{Cata};
recall that, by definition, hyperelliptic and Enriques
surfaces are respectively quotients of Abelian and $K3$-surfaces by
free involutions). Furthermore, quasi-simplicity of hyperelliptic and
Enriques surfaces extends to quasi-simplicity of the quotients of Abelian
and $K3$-surfaces by certain finite group actions, see \cite{Duke}.

Whether elliptic surfaces and irrational ruled non $\mathbb
C$-minimal
surfaces are quasi-simple is, as far as I know, an open question.

It is natural to expect that for surfaces of general type there is no
quasi-simplicity: there should exist examples of real deformation
distinct real surfaces with diffeomorphic real structures.
A challenging problem is to find convenient deformation invariants
which are not covered by the differential topology of $(X_\C, c)$.

Existence of non quasi-simple families of surfaces
of general type does not prevent certain particular classes of surfaces
of general type from being quasi-simple. And examples of quasi-simple
real surfaces of general type do exist. One such example is given by real
Bogomolov-Miyaoka-Yau surfaces, that is, surfaces covered by a ball in
$\C^2$, see \cite{KK}. (Note in passing that in \cite{KK} it is also shown
that there exist diffeomorphic, in fact complex conjugated,
Bogomolov-Miyaoka-Yau surfaces which are not real and thus, being
rigid, are not deformation equivalent as complex surfaces.
These surfaces are counter-examples
to the so called $\rm Diff=Deff$ problem in complex geometry, see
\cite{KK} for precise
definitions and references to counter-examples not related to
the complex conjugation.
This problem is a kind of substitute of
quasi-simplicity for complex varieties. The existence of
$\rm Diff\ne Deff$ examples explains why we need to fix complex deformation
class in the definition of quasi-simplicity of real
varieties.)\footnote{{\it Added in proof.} When this paper had
been
already finished, we with
Vik.~Kulikov have constructed examples of non quasi-simple real surfaces of general type.}

\subsection{ Finiteness}

While the problem of quasi-simplicity is solved for rational
surfaces and is essentially open and very difficult for surfaces of
general type, the situation with finiteness is an opposite one:
finiteness holds both for each complex deformation class of surfaces
of general type (deformation finiteness)
and for any fixed surface of general type (individual finiteness).

To wit, since the composition of two real structures is a biholomorphic
automorphism and since the group of automorphisms of any variety of general
type is finite, there are only finitely many real structures on a
variety of general type (the same argument works for nonsingular hypersurfaces
of degree $\ge 3$ in projective space of dimension $n\ge 3$ with the exception
of $n=3, d=4$).
This is what we call {\it individual finiteness}, which we understand
as finiteness of the number of conjugacy classes of real structures on a given variety
(note that individual finiteness understood in this way extends to hypersurfaces
of degree $4$ in projective spaces of dimension $3$, see \cite{DIK}).

On the other hand, the Hilbert scheme of varieties of general type with given
characteristic numbers is quasi-projective, which implies {\it deformation
finiteness:} real structures on the varieties which, as complex varieties, are
deformation equivalent to a given variety of general type split into a finite number
of real deformation classes (where, according to our definitions, both variety
and real structure are subject to deformation).

Unlike surfaces of general type, a rational surface may have a huge automorphism
group, and, as far as I know, the problem of individual finiteness for rational
surfaces is open. The situation is different with regard to deformation finiteness
of rational surfaces which is an easy consequence of their quasi-simplicity.

In fact, {\it the deformation finiteness holds for any type of
surfaces}. Indeed, the only birational classes of surfaces for which such a
result is not contained in the literature, either explicitly or implicitly,
are elliptic surfaces and irrational ruled surfaces, but for these classes the
proof is more or less straightforward. It would be useful to find a conceptual
proof dealing with all types of surfaces in a unified way.

Some finiteness results are also known
for Klein actions of finite
groups on $K3$- and Abelian surfaces.
In particular,
{\it the number of equivariant deformation classes of $K3$- and Abelian surfaces with
faithful Klein actions of finite groups is finite}, see~\cite{Duke}.

Another, higher-dimensional, generalization of
finiteness of real structures on $K3$-surfaces extends it to
so called holomorphic symplectic (hence hyperk\"ahler) manifolds:
{\it the number of equivariant deformation classes of real structures in a given
deformation class of compact holomorphic symplectic manifolds is finite},
see \cite{inpreparation}.)

The differential topology of $(X_\C, c)$ is preserved under
deformation, and therefore deformation finiteness implies
topological finiteness. Another, more direct, approach to topological
finiteness was  recently developed by Y.~Laszlo and C.~Viterbo.
They proposed to study  finiteness of diffeomorphism types
of  real forms on complex projective varieties of a given degree.
Here one should distinguish between the real and complex degree.
For example, there exists a sequence $X_n$ of complex $K3$-surfaces
of degree $4$ (quartics in $\PP^3$) such that, for appropriate
real structures $c_n$ on $X_n$, their real degrees (the minimal
degree of a real projective embedding $X_n\to\PP^{q_n}$) converge to infinity
(so that these real structures are not induced from $\PP^3$ and, moreover,
can not be induced from $\PP^q$ with bounded $q$).
Of course, varieties of a given real degree split into a finite
number of families. Whether the same is true for real varieties of a given
complex degree is still an open question of utmost importance. But, thanks to
Y.~Laszlo and C.~Viterbo \cite{LV}, we  now have some explicit bounds
for the Pontryagin numbers of varieties of a given real degree and, as a
consequence, some explicit bounds for the number of cobordism classes of such
varieties.


\begin{thebibliography}{99}

\bibitem{Cata} Catanese, F., and  Frediani, P.:
Real hyperelliptic surfaces and the orbifold fundamental group. {\it J. Inst. Math. Jussieu}, {\bf 2}
(2003), 163--223.




\bibitem{DIK} Degtyarev, A., Itenberg, I., and  Kharlamov, V.:
Real Enriques Surfaces. {\it Lecture Notes in Math., Springer, Berlin}, {\bf 1746}
(2000).

\bibitem{Duke} Degtyarev, A., Itenberg, I., and  Kharlamov, V.:
Finiteness and quasi-simplicity for symmetric $K3$-surfaces. {\it Duke Math. J.}, {\bf 122}
(2004), no. 1, 1--49.

\bibitem{inpreparation} Degtyarev, A., Itenberg, I., and  Kharlamov, V.:
Finiteness for real hyperk\"ahler manifolds. {\it in preparation}.



\bibitem{DK} Degtyarev, A., and  Kharlamov, V.:
Topological properties of real algebraic
varieties: Rokhlin's way. {\it Russ. Math. Surveys}, {\bf 55}
(2000), no. 4, 735--814.

\bibitem{DKrat} Degtyarev, A., and  Kharlamov, V.:
Real rational surfaces are quasi-simple. {\it J. Reine. Angew. Math.}, {\bf 551}
(2002), 87--99.

\bibitem{Kh} Kharlamov, V.:
Generalized Petrovskii inequality.
{\it Funkz. Anal. i Priloz.},
{\bf 9} (1974), p.  50--56.


\bibitem{Kh2} Kharlamov, V.:
Generalized Petrovskii inequality II.
{\it Funkz. Anal. i Priloz.},
{\bf 10}
(1975),
p.  93--94.

\bibitem{Kh3} Kharlamov, V.:
Topology of real algebraic varieties.
{\it in Collected Papers by Petrovskii, Nauka}, 1986,
p. 546--598.

\bibitem{Khchir} Kharlamov, V.:
On non-amphicheiral surfaces of degree 4 in ${\Bbb R} P^3$.
{\it Lecture Notes in Math.}
{\bf 1346}
(1988), p. 349--356.




\bibitem{Bour} Kharlamov, V.:
Vari\'et\'es de Fano r\'eelles (d'apr\`es C.~Viterbo). {\it Ast\'erisque, S\'eminaire Bourbaki}, {\bf 276}
(2002), p. 189--206.

\bibitem{KK} Kharlamov, V., and Kulikov, Vik.:
On real structures of rigid surfaces.
{\it Izv. Math.}, {\bf 66} (2003), no. 1,
p. 133--150.

\bibitem{Ko1} Koll\'ar J.:
The topology of real and complex algebraic varieties.
{\it Mathematical Society of Japan. Adv. Stud. Pure Math.} {\bf 31} (2001),
p. 127--145

\bibitem{Ko2} Koll\'ar J.:
The Nash conjecture for nonprojective threefolds.
{\it Contemp. Math.} {\bf 312} (2002),
p. 137-152.

\bibitem{LV} Laszlo, Y.
and Viterbo, C.: Estimates of characteristic numbers of real algebraic
varieties. {\it preprint} (2004),
1--24.

\bibitem{Man} Manin, Yu.I. and Tsfasman, M.A.:
Rational varieties: Algebra, geometry and arithmetic.
{\it Russ. Math. Surv.}, {\bf 41} (1986), no. 2,
p. 51--116.

\bibitem{Mi} Milnor, J.:
On the Betti numbers of real
varieties. {\it Proc. Amer. Math. Soc.} {\bf 15} (1964),
275--280.

\bibitem{Nik} Nikulin, V.:
Kummer surfaces.
{\it Math. USSR - Izv.}, {\bf 9} (1975), no. 2,
p. 261--275.


\bibitem{Nik2} Nikulin, V.:
Integer symmetric bilinear forms and some of their geometric applications.
{\it Math. USSR - Izv.}, {\bf 14} (1979), no. 1,
p. 103--167.

\bibitem{Th} Thom, R.:
Sur l'homologie des vari\'et\'es alg\'ebriques r\'eelles. {\it in
Differential and Combinatorial Topology. Symp. Marston Morse} (1965),
255--265.

\bibitem{Var} Varchenko, A.N.:
On a local residue and the intersection form in vanishing cohomologies.
{\it in Izv. Akad. Nauk SSSR, Ser. Mat.}
{\bf  49},
(1985),
p. 32--54.



\bibitem{W} Welschinger, J.Y.:
Real structures on minimal ruled sufaces.
{\it Comment. Math. Helv.}, {\bf 78} (2003),
p. 418--466.

%
%
%
%
%

\end{thebibliography}
\end{document}